\documentclass[10pt,twoside,reqno]{amsart}
\usepackage{amssymb}
\usepackage{multirow}
\calclayout

\numberwithin{equation}{section}
\makeatletter

\renewcommand{\@secnumfont}{\bfseries}

\renewcommand{\section}{\@startsection{section}{1}%
  {0mm}{.7\linespacing\@plus\linespacing}{.5\linespacing}
  {\normalfont\bfseries\centering}}

\newcommand{\bibsection}{\@startsection{section}{1}%
  {0mm}{.7\linespacing\@plus\linespacing}{.5\linespacing}
  {\normalfont\scshape\centering}}

\renewcommand{\@biblabel}[1]{#1.}

\begin{document}

\vspace{1.3cm}

\title {A note on Catalan numbers associated with $p$-adic integral on $\mathbb{Z}_p$}

\author{Taekyun Kim}
\address{Department of Mathematics \\ Kwangwoon University \\ Seoul 139-701 \\ Republic
	of Korea}
\email{tkkim@kw.ac.kr}

\thanks{\scriptsize }

\maketitle

\begin{abstract}
In this paper, we study Catalan numbers which can be represented by the $p$-adic integral on $\mathbb{Z}_p$ and we investigate some properties and formulae related to Catalan numbers and special numbers.
\end{abstract}

\bigskip
\medskip
\section{\bf Introduction}
Let $p$ be a f
ixed odd prime number. Throughout this paper, $\mathbb{Z}_p$,$\mathbb{Q}_p$ and $\mathbb{C}_p$ will denote the ring of $p$-adic integers, the field of $p$-adic numbers and the completion of the algebraic closure of $\mathbb{Q}_p$. The $p$-adic norm $|\cdot|_p$ is normalized as $|p|_p = \frac{1}{p}$. As is well known, the Euler numbers are defined by the generating function to be
\begin{equation}\begin{split}\label{01}
\frac{2}{e^t+1}e^{xt} = \sum_{n=0}^\infty E_n(x) \frac{t^n}{n!}, \quad (\textnormal{see} \,\, [1-13]).
\end{split}\end{equation}
When $x=0$, $E_n=E_n(0)$ are called the Euler numbers. From \eqref{01}, we note that
\begin{equation}\begin{split}\label{02}
E_n(x) = \sum_{l=0}^n {n \choose l} x^{n-l} E_l.
\end{split}\end{equation}
Recently, $\lambda $-Changhee polynomials are defined by the generating function to be
\begin{equation}\begin{split}\label{03}
\frac{2}{(1+t)^{\lambda }+1}(1+t)^{\lambda x} = \sum_{n=0}^\infty Ch_{n,\lambda }(x) \frac{t^n}{n!},\quad\textnormal{where}\,\, \lambda  \in \mathbb{Z}_p \quad (\textnormal{see}\,\,[5]).
\end{split}\end{equation}
By replacing $t$ by $e^{\frac{1}{\lambda }t}-1$, we get
\begin{equation}\begin{split}\label{04}
\frac{2}{e^t+1}e^{xt} =& \sum_{n=0}^\infty Ch_{n,\lambda }(x) \frac{1}{n!} \Big( e^{\frac{1}{\lambda }t}-1\Big)^n \\
=& \sum_{n=0}^\infty Ch_{n,\lambda }(x) \sum_{m=n}^\infty S_2(m,n) \lambda ^{-m} \frac{t^m}{m!},
\end{split}\end{equation}
where $S_2(m,n)$ is the stirling number of the second kind. When $x=0$, $Ch_{n,\lambda } = Ch_{n,\lambda }(0)$ are called the $\lambda $-Changhee numbers. From \eqref{04}, we note that
\begin{equation}\begin{split}\label{05}
E_m(x) = \sum_{n=0}^m Ch_{n,\lambda }(x) S_2(m,n) \lambda ^{-m}, \quad(m \geq 0).
\end{split}\end{equation}
For $n \geq 0$, the stirling number of the first kind is defined as
\begin{equation}\begin{split}\label{06}
(x)_n = x(x-1)\cdots(x-n+1) = \prod_{l=1}^{n-1} (x-l) = \sum_{l=0}^n S_1(n,l) x^l,
\end{split}\end{equation}
and the stirling number of the second kind is given by
\begin{equation}\begin{split}\label{07}
x^n = \sum_{l=0}^n S_2(n,l) (x)_l, \quad (\textnormal{see} \,\, [7,8,9]).
\end{split}\end{equation}
Let $C(\mathbb{Z}_p)$ be the space of all continuous functions on $\mathbb{Z}_p$. For $f \in C(\mathbb{Z}_p)$, the fermionic $p$-adic integral on $\mathbb{Z}_p$ is defined by Kim to be
\begin{equation}\begin{split}\label{08}
\int_{\mathbb{Z}_p} f(x) d \mu_{-1}(x) =& \lim_{N \rightarrow \infty} \sum_{x=0}^{p^N-1} f(x) \mu_{-1} (x+ p^N \mathbb{Z}_p)\\
=& \lim_{N \rightarrow \infty} \sum_{x=0}^{p^N-1} f(x) (-1)^x, \quad (\textnormal{see} \,\, [8]).
\end{split}\end{equation}
From \eqref{08}, we note that
\begin{equation}\begin{split}\label{09}
I_{-1}(f_1) + I_{-1}(f) = 2f(0), \quad \textnormal{where}\,\, f_1(x) = f(x+1).
\end{split}\end{equation}
Thus, by \eqref{09}, we get
\begin{equation}\begin{split}\label{10}
I_{-1}(f_n) + (-1)^{n-1}I_{-1}(f) = 2 \sum_{l=0}^{n-1}(-1)^{n-1-l}f(l), \quad (\textnormal{see} \,\, [8]).
\end{split}\end{equation}
where $f_n(x) = f(x+n),\,\,n \in \mathbb{N}$. It is not difficult to show that
\begin{equation}\begin{split}\label{11}
\frac{1-\sqrt{1-4t}}{2t}=& \frac{(-1)}{2t} \left( \sum_{n=1}^\infty {2n \choose n} \frac{(-1)^{n-1}}{4^n(2n-1)}(-1)^n 4^n t^n \right)       \\
=& \frac{1}{2} \sum_{n=1}^\infty {2n \choose n} \frac{1}{2n-1}t^{n-1}       \\
=& \frac{1}{2} \sum_{n=0}^\infty {2n+2 \choose n+1} \frac{1}{2n+1} t^n       \\
=& \frac{1}{2} \sum_{n=0}^\infty \frac{(2n+2)(2n+1)(2n)\cdots(n+2)}{(n+1)n!(2n+1)}t^n              \\
=& \sum_{n=0}^\infty {2n \choose n} \frac{1}{n+1} t^n,\quad (\textnormal{see} \,\, [1,2,3]).      \\
\end{split}\end{equation}
As is well known, the Catalan number $C_n$ is defined by $C_n = {2n \choose n} \frac{1}{n+1}, \,\,(n \geq 0).$ From \eqref{11}, we note that the generating function of Catalan numbers is given by
\begin{equation}\begin{split}\label{12}
\frac{2}{1+\sqrt{1-4t}} = \sum_{n=0}^\infty n! C_n \frac{t^n}{n!},\quad (\textnormal{see} \,\, [9-13]).
\end{split}\end{equation}
In this paper, we study Catalan numbers associated with $p$-adic integral on $\mathbb{Z}_p$ and we give Witt's type formula related Catalan numbers.
\section{Catalan numbers associated with $p$-adic intergral on $\mathbb{Z}_p$}
For $t \in \mathbb{C}_p$ with $|t|_p < p^{-\frac{1}{p-1}}$, we observe that
\begin{equation}\begin{split}\label{13}
\int_{\mathbb{Z}_p} (1+t)^{\tfrac{x}{2}}   d\mu_{-1} (x) = \frac{2}{1+\sqrt{1+t}} = \sum_{n=0}^\infty Ch_{n,\frac{1}{2}} \frac{t^n}{n!}.
\end{split}\end{equation}
From \eqref{09}, we note that
\begin{equation}\begin{split}\label{14}
\int_{\mathbb{Z}_p}  e^{xt}  d\mu_{-1} (x) = \frac{2}{e^t+1} \sum_{n=0}^\infty E_n \frac{t^n}{n!}
\end{split}\end{equation}
On the other hand,
\begin{equation}\begin{split}\label{15}
\int_{\mathbb{Z}_p} (1+t)^{\tfrac{x}{2}}   d\mu_{-1} (x) &=  \sum_{m=0}^\infty \int_{\mathbb{Z}_p} x^m   d\mu_{-1} (x) \frac{1}{2^m} \frac{1}{m!} \Big(\log(1+t)\Big)^m
\\&=\sum_{m=0}^\infty E_m 2^{-m} \sum_{n=m}^\infty S_1(n,m) \frac{t^n}{n!}
\\&= \sum_{n=0}^\infty \left( \sum_{m=0}^n 2^{-m} E_m S_1(n,m) \right) \frac{t^n}{n!}.
\end{split}\end{equation}
Therefore, by \eqref{13} and \eqref{15}, we obtain the following theorem. \\\\
\textbf{Theorem 1}. For $n \geq 0$, we have
\begin{equation*}
\sum_{m=0}^n 2^{-m} E_m S_1(n,m) = Ch_{n, \frac{1}{2}}
\end{equation*}
By replacing $t$ by $-4t$ in \eqref{13}, we get
\begin{equation}\begin{split}\label{16}
\int_{\mathbb{Z}_p}  (1-4t)^{\tfrac{x}{2}}  d\mu_{-1} (x) =& \frac{2}{1+\sqrt{1-4t}} = \sum_{n=0}^\infty \frac{1}{n+1} {2n \choose n} t^n \\
=& \sum_{n=0}^\infty C_n t^n.
\end{split}\end{equation}
On the other hand,
\begin{equation}\begin{split}\label{17}
\int_{\mathbb{Z}_p}  (1-4t)^{\tfrac{x}{2}}  d\mu_{-1} (x) = \sum_{n=0}^\infty Ch_{n,\frac{1}{2}} \frac{(-4)^n t^n}{n!}.
\end{split}\end{equation}
Therefore, by \eqref{16} and \eqref{17}, we obtain the following theorem.\\\\
\textbf{Theorem 2.} For $n \geq 0$, we have
\begin{equation*}
C_n = (-1)^n \frac{4^n Ch_{n,\frac{1}{2}}}{n!}.
\end{equation*}
From \eqref{16}, we have
\begin{equation}\begin{split}\label{18}
\sum_{n=0}^\infty (-1)^n 4^n \int_{\mathbb{Z}_p} {\frac{x}{2} \choose n}   d\mu_{-1} (x)t^n &= \sum_{n=0}^\infty C_n t^n \\
&= \sum_{n=0}^\infty \frac{1}{n+1} {2n \choose n} t^n.
\end{split}\end{equation}
Therefore, by comparing the coefficients on the both sides of \eqref{18}, we obtain the following theorem.\\\\
\textbf{Theorem 3.} For $n \geq 0$, we have
\begin{equation*}
\int_{\mathbb{Z}_p} {\frac{x}{2} \choose n}    d\mu_{-1} (x)= \frac{(-1)^n}{4^n} C_n = \frac{(-1)^n}{4^n (n+1)} {2n \choose n}.
\end{equation*}
From \eqref{15}, we have the following equation:
\begin{equation}\begin{split}\label{19}
\int_{\mathbb{Z}_p} {\frac{x}{2} \choose n}   d\mu_{-1} (x) = \frac{1}{n!} \sum_{m=0}^n 2^{-m} E_m S_1(n,m)
\end{split}\end{equation}
Therefore, by Theorem 3 and \eqref{19}, we obtain the following theorem.\\\\
\textbf{Theorem 4.} For $n \geq 0$, we have
\begin{equation*}
C_n = \frac{(-1)^n}{n!} \sum_{m=0}^n 2^{2n-m} E_m S_1(n,m).
\end{equation*}
Now, we observe that
\begin{equation}\begin{split}\label{20}
\sqrt{1+t}&= (1+t)^{\tfrac{1}{2}} = \sum_{n=0}^\infty {\frac{1}{2} \choose n} t^n = \sum_{n=0}^\infty \frac{\big( \frac{1}{2}\big)_n}{n!} t^n\\
&= \sum_{n=0}^\infty  \frac{1 \big(\tfrac{1}{2}-1 \big) \big( \tfrac{1}{2}-2\big) \cdots \big( \tfrac{1}{2}-n+1 \big) }{n!} t^n \\
&=\sum_{n=0}^\infty	 \frac{(-1)^{n-1} 1 \cdot 3 \cdot 5 \cdots (2n-3)}{n! 2^n} t^n\\
&= \sum_{n=0}^\infty \frac{(-1)^{n-1} 1 \cdot 2 \cdot 3 \cdot 4 \cdots (2n-3) (2n-2) (2n-1)(2n)}{n! 2^n 2 \cdot 4 \cdot 6 \cdots (2n-2) (2n-1) (2n)} t^n\\
&= \sum_{n=0}^\infty (-1)^{n-1} \frac{(2n)!}{n! 4^n (2n-1)n!} t^n = \sum_{n=0}^\infty {2n \choose n} \frac{(-1)^{n-1}}{4^n (2n-1)} t^n.
\end{split}\end{equation}
By \eqref{13} and \eqref{20}, we get
\begin{equation}\begin{split}\label{21}
2 &= \left( \sum_{n=0}^\infty Ch_{n,\frac{1}{2}} \frac{t^n}{n!} \right) \Big( 1+ \sqrt{1+t} \Big) \\
&= \sum_{n=0}^\infty Ch_{n,\frac{1}{2}} \frac{t^n}{n!} + \left( \sum_{k=0}^\infty Ch_{k,\frac{1}{2}} \frac{t^k}{k!} \right) \left( \sum_{m=0}^\infty {2m \choose m} \frac{(-1)^{m-1}}{4^m (2m-1)} t^m \right)\\
&= \sum_{n=0}^\infty Ch_{n,\frac{1}{2}} \frac{t^n}{n!} + \sum_{n=0}^\infty \left( \sum_{m=0}^n C_m \frac{(m+1)(-1)^{m-1}}{4^m(2m-1)} \frac{m!n!}{(n-m)!m!} Ch_{n-m, \frac{1}{2}} \right) \frac{t^n}{n!}\\
&= \sum_{n=0}^\infty Ch_{n,\frac{1}{2}} \frac{t^n}{n!} + \sum_{n=0}^\infty \left( \sum_{m=0}^n C_m \frac{(m+1)!(-1)^{m+1}}{4^m(2m-1)}{n \choose m}Ch_{n-m, \frac{1}{2}} \right) \frac{t^n}{n!}\\
\end{split}\end{equation}
By comparing the coefficients on the both sides of \eqref{21}, we get
\begin{equation}\begin{split}\label{22}
Ch_{n,\frac{1}{2}} + \sum_{m=0}^n {n \choose m} C_m Ch_{n-m, \frac{1}{2}} (m+1)!(-1)^{m+1}
\frac{1}{4^m(2m-1)} = \begin{cases}
2 \quad \textnormal{if} \,\,n=0
\\
0 \quad \textnormal{if} \,\,n >0.
\end{cases}
\end{split}\end{equation}
Therefore, by \eqref{22}, we obtain the following theorem.\\\\
\textbf{Theorem 5.} For $n \in \mathbb{N}$, we have
\begin{align*}
&Ch_{n,\frac{1}{2}} = \sum_{m=0}^n {n \choose m} C_m Ch_{n-m, \frac{1}{2}}(-1)^{m}
\frac{(m+1)! }{4^m(2m-1)}, \\
\textnormal{and}& \\
&Ch_{0,\frac{1}{2}} = 1.
\end{align*}
By replacing $t$ by $\tfrac{1}{4} \big( 1-e^{2t} \big)$ in \eqref{12}, we get
\begin{equation}\begin{split}\label{23}
\frac{2}{e^t+1} &= \sum_{n=0}^\infty C_n \frac{1}{4^n} \big( 1- e^{2t} \big)^n \\
&= \sum_{n=0}^\infty C_n \frac{(-1)^n}{4^n} \big( e^{2t}-1 \big)^n \\
&= \sum_{n=0}^\infty C_n \frac{(-1)^n}{4^n} n! \sum_{m=n}^\infty S_2(m,n) \frac{2^m t^m}{m!}\\
&= \sum_{m=0}^\infty \left( \sum_{n=0}^m C_n (-1)^n 2^{m-2n} n! S_2(m,n) \right) \frac{t^m}{m!}.
\end{split}\end{equation}
Therefore, by \eqref{01} and \eqref{23}, we obtain the following theorem.\\\\
\textbf{Theorem 6.} For $m \geq 0$, we have
\begin{equation*}
E_m = \sum_{n=0}^m C_n (-1)^n 2^{m-2n} n! S_2(m,n).
\end{equation*}
Now, we observe
\begin{equation}\begin{split}\label{24}
(1+t)^{\tfrac{x}{2}} &= \sum_{m=0}^\infty \Big( \frac{x}{2} \Big)^m \frac{\big(\log(1+t)\big)^m}{m!} = \sum_{m=0}^\infty  \Big( \frac{x}{2} \Big)^m \sum_{n=m}^\infty S_1(n,m) \frac{t^n}{n!} \\
&= \sum_{n=0}^\infty \left( \sum_{m=0}^n \Big( \frac{x}{2} \Big)^m S_1(n,m) \right) \frac{t^n}{n!}.
\end{split}\end{equation}
We consider the $\frac{1}{2}$-Changhee polynomials which are given by the generating function to be
\begin{equation}\begin{split}\label{25}
\int_{\mathbb{Z}_p} (1+t)^{\tfrac{x+y}{2}}   d\mu_{-1} (y) &= \frac{2}{1+\sqrt{1+t}} \sqrt{(1+t)^x}\\
&= \sum_{n=0}^\infty Ch_{n,\frac{1}{2}} (x) \frac{t^n}{n!}.
\end{split}\end{equation}
When $x=0$, we note that $Ch_{n,\frac{1}{2}} = Ch_{n,\frac{1}{2}}(0)$.
By \eqref{23}, \eqref{24} and \eqref{25}, we  get
\begin{equation}\begin{split}\label{26}
\sum_{n=0}^\infty Ch_{n,\frac{1}{2}} (x) \frac{t^n}{n!} &= \left( \frac{2}{1+\sqrt{1+t}} \right) \left( (1+t)^{\frac{x}{2}} \right)\\
&= \left( \sum_{k=0}^\infty Ch_{k,\frac{1}{2}} \frac{t^k}{k!} \right) \left( \sum_{m=0}^\infty \Bigg( \sum_{j=0}^m \Big(\frac{x}{2} \Big)^j S_1(m,j) \Bigg) \frac{t^m}{m!} \right)\\
&= \sum_{n=0}^\infty \left( \sum_{m=0}^n \sum_{j=0}^m \Big( \frac{x}{2} \Big)^j S_1(m,j) Ch_{n-m, \frac{1}{2}} {n \choose m} \right) \frac{t^n}{n!}.
\end{split}\end{equation}
Therefore, by \eqref{26}, we obtain the following theorem.\\\\
\textbf{Theorem 7.} For $n \geq 0$, we have
\begin{equation*}
Ch_{n,\frac{1}{2}}(x) =  \sum_{m=0}^n \sum_{j=0}^m \Big( \frac{x}{2} \Big)^j S_1(m,j) Ch_{n-m, \frac{1}{2}} {n \choose m}.
\end{equation*}
By replacing $t$ by $-4t$ in \eqref{25}, we define Catalan polynomials which are given by the generating function to be
\begin{equation}\begin{split}\label{27}
\int_{\mathbb{Z}_p}  (1-4t)^{\tfrac{x+y}{2}}  d\mu_{-1} (y) &= \frac{2}{1+\sqrt{1-4t}} \sqrt{(1-4t)^x} \\
&= \sum_{n=0}^\infty C_n(x) t^n.
\end{split}\end{equation}
From \eqref{12} and \eqref{24}, we note that
\begin{equation}\begin{split}\label{28}
\sum_{n=0}^\infty C_n(x) t^n &= \left( \frac{2}{1+\sqrt{1-4t}} \right) \left( (1-4t)^{\tfrac{x}{2}} \right) \\
&= \left( \sum_{k=0}^\infty C_k t^k \right) \left( \sum_{m=0}^\infty \Bigg( \sum_{j=0}^m \Big( \frac{x}{2} \Big)^j S_1(m,j) \Bigg) \frac{(-4)^m t^m}{m!}  \right)\\
&= \sum_{n=0}^\infty \left( \sum_{m=0}^n \sum_{j=0}^m \Big( \frac{x}{2} \Big)^j S_1(m,j) \frac{(-4)^m}{m!} C_{n-m} \right) t^n.
\end{split}\end{equation}
Therefore, by comparing the coefficients on the both sides of \eqref{28}, we obtain the following theorem. \\\\
\textbf{Theorem 8.} For $n \geq 0$, we have
\begin{equation*}
C_n(x) = \sum_{m=0}^n \sum_{j=0}^m \Big( \frac{x}{2} \Big)^j S_1(m,j) \frac{(-4)^m}{m!} C_{n-m} .
\end{equation*}
Remark. From \eqref{20} and \eqref{24}, we note that
\begin{equation}\begin{split}\label{29}
\sum_{n=0}^\infty {2n \choose n} \frac{(-1)^{n-1}}{4^n (2n-1)} t^n = \sum_{n=0}^\infty \left( \frac{1}{n!} \sum_{m=0}^n \Big( \frac{1}{2} \Big)^m S_1(n,m) \right) t^n
\end{split}\end{equation}
By \eqref{29}, we get
\begin{equation}\begin{split}\label{30}
\frac{1}{n+1} {2n \choose n} \frac{(-1)^{n-1}n+1 }{4^n (2n-1)} = \frac{1}{n!} \sum_{m=0}^n \Big(\frac{1}{2}\Big)^m S_1(n,m)
\end{split}\end{equation}
Thus, from \eqref{30}, we have
\begin{equation*}
C_n = \frac{4^n(2n-1)}{(n+1)!} (-1)^{n-1} \sum_{m=0}^n \Big( \frac{1}{2} \Big)^m S_1(n,m), \,\, (n \geq 0).
\end{equation*}
\textbf{Corollary 9.} For $n \geq 0$, we have
\begin{equation*}
C_n = (-1)^{n+1} \frac{4^n (2n-1)}{(n+1)!} \sum_{m=0}^n \Big( \frac{1}{2} \Big)^m S_1(n,m).
\end{equation*}

\end{document}